\documentclass[12pt]{article}
\usepackage{amsmath}
\usepackage{amstext}
\usepackage{amsfonts}
\usepackage{amsthm}
\usepackage{amscd}
\usepackage{amssymb}
\numberwithin{equation}{section}

\swapnumbers
\theoremstyle{plain}
\newtheorem{thm}{Theorem}[section]
\newtheorem{lem}[thm]{Lemma}
\newtheorem{prop}[thm]{Proposition}
\newtheorem{cor}[thm]{Corollary}
\theoremstyle{remark}

\def\mod0{{\mathcal SU}_C(2,{\mathcal O})}

\def\modp{{\mathcal SU}_C(2,{\mathcal O}(p))}

\def\cO{{\mathcal O}}

\def\map#1{\ \smash{\mathop{\longrightarrow}\limits^{#1}}\ }

\def\LL{\mathcal{L}}

\def\tC{\tilde{C}}
\def\tp{\tilde{p}}

\def\tk{\tilde{\kappa}}

\def\sym{\mathrm{Sym}}

\def\pic{\mathrm{Pic}}

\def\im{\mathrm{im} \:}
\def\inj{\hookrightarrow}
\def\Spin{\mathrm{Spin}}
\def\SO{\mathrm{SO}}
\def\SL{\mathrm{SL}}
\def\GL{\mathrm{GL}}
\def\Hom{\mathrm{Hom}}
\def\Nm{\mathrm{Nm}}
\def\id{\mathrm{id}}

\def\lra{\longrightarrow}

\def\ra{\rightarrow}

\def\lms{\longmapsto}

\def\cc{\mathbb{C}}
\def\pp{\mathbb{P}}

\def\MM{\mathcal{M}}
\def\DD{\mathcal{D}}
\def\HH{\mathcal{H}}
\def\tC{\tilde{C}}
\def\EE{\mathcal{E}}
\def\DET{\mathrm{DET}}
\def\PF{\mathrm{PF}}
\def\div{\mathrm{div}}
\def\Heis{\mathrm{Heis}}
\def\End{\mathrm{End}}
\def\supp{\mathrm{supp}}

\title{A duality for Spin Verlinde spaces and Prym theta
functions}
\author{C. Pauly and S. Ramanan}

\begin{document}

\maketitle

\begin{abstract}
We prove canonical isomorphisms between Spin Verlinde spaces,i.e.,
spaces of global sections of a determinant line bundle over the
moduli space of semistable $\Spin_n$-bundles over a smooth
projective curve $C$, and the dual spaces of theta functions
over Prym varieties of unramified double covers of $C$.
\end{abstract}

\section{Introduction}
To any smooth, projective  curve $C$, one classically 
associates a collection
of principally polarized abelian varieties: the Jacobian $JC$,
parametrizing degree zero line bundles, as well as, for
any unramified double cover $\tC_\eta \ra C$, depending on a
$2$-torsion point $\eta \in JC[2]$, a Prym variety $P_\eta$. The
projective geometry of the configuration $JC \cup \bigcup_\eta P_\eta$ 
has been much studied \cite{mum1}, \cite{VGP1}, \cite{beau1} and
encodes e.g. the Schottky-Jung identities among theta-constants 
\cite{mum1}.

Less classically, one can consider the moduli space $\MM(G)$ of
semistable principal $G$-bundles over the curve $C$, where
$G$ is a simple and simply-connected algebraic group. For some
ample line bundle $\LL$ over $\MM(G)$, the vector space of global 
sections $H^0(\MM(G),\LL)$ has been identified to a space
of conformal blocks arising in conformal field theory 
(see e.g.\cite{sorger} for a survey or \cite{ls} for a proof), which
made the computation of its dimension possible. This is the
celebrated Verlinde formula.

In this paper we are interested in dualities between Verlinde
spaces $H^0(\MM(G),\LL)$ and spaces of abelian theta functions
,i.e.,sections of some multiple of a principal polarization over $JC$ or,
more generally, $JC \cup \bigcup_\eta P_\eta$. Such dualities were
first proved for the structure group $G = \SL_2$ and the line
bundles $\LL^l$ for $l=1,2,4$ \cite{beau1}, \cite{beauville}, 
\cite{op} or, more generally, $G = \SL_n$ and the line
bundle $\LL$ \cite{bnr}, where $\LL$ is the ample generator of the
Picard group of $\MM(\SL_n)$.

In the articles \cite{oxbury1} and \cite{oxbury2}, W.M. Oxbury 
constructs linear maps between a Verlinde space for the complex
spin group,i.e.,$G = \Spin_n$, and a space of abelian theta
functions over $JC \cup \bigcup_\eta P_\eta$. Our main theorem
states that these linear maps are actually isomorphisms. More
precisely we show

\begin{thm} \label{duality}
For any curve $C$ and any integer $m\geq 1$, we have canonical 
isomorphisms
\begin{equation} \label{prymspinodd}
(1) \sum s_\eta^\pm \ : \  \sum_{\eta \in JC[2]} 
 H^0(P^{ev}_\eta, (2m+1)\Xi_\eta)^\vee_\pm
\map{\sim} H^0(\MM^\pm(\Spin_{2m+1}), \Theta(\cc^{2m+1}))
\end{equation}
\begin{equation} \label{prymspineven}
(2) \sum s_\eta^\pm \ : \  \sum_{\eta \in JC[2]}
\left[ H^0(P^{ev}_\eta, (2m)\Xi_\eta)^\vee_\pm
\oplus H^0(P^{odd}_\eta, (2m)\Xi_\eta)^\vee_\pm \right]
\map{\sim} H^0(\MM^\pm(\Spin_{2m}), \Theta(\cc^{2m}))
\end{equation}
\end{thm}

We refer to section 2 and 3 for the rather technical details encoded
in the notation.

We note that theorem 1.1 for the group $\Spin_3$ coincides with the
above mentioned duality for the Verlinde space $H^0(\MM(\SL_2),
\LL^4)$ under the
exceptional isomorphism $\Spin_3 = \SL_2$. In this case theorem
1.1 was proved in \cite{ram2} for all curves and in \cite{op}
for the generic curve.

The proof of the theorem essentially goes as follows: first,
we construct maps from products of Prym varieties, called $SCH_n$, into
the moduli space $\MM(\Spin_n)$. This is done by induction on $n$,
exploiting the exceptional isomorphisms of $\Spin_n$ with algebraic
groups for small values of $n$. For example, $SCH_3$ is isomorphic
to the union $JC \cup \bigcup_\eta P_\eta$. Next, we observe that the
divisor in the product 
$\MM(\Spin_n) \times (JC \cup \bigcup_\eta P_\eta)$, which
induces the pairing of the theorem, decomposes when restricted
to the variety $SCH_n$. Hence we obtain a factorization of the
linear map of theorem 1.1. As our auxiliary variety $SCH_n$ is
constructed from abelian varieties, we can use theorems on 
multiplication of theta functions (proposition 4.1 and 4.2) 
and dualities for
second order theta functions (Prym-Wirtinger duality, proposition 5.1) 
in order to show, again inductively, that all maps of the
factorization are injective. By the Verlinde formula,
this will be enough to deduce the theorem.

It is legitimate to ask whether the isomorphism of theorem 1.1 is
a particular case of more general duality. Indeed, in the case
$G = \SL_n$, one constructs (\cite{beau3} section 8 and \cite{dt}),
via canonical theta divisors,
pairings among $\SL_n$ (and $\GL_n$)-Verlinde spaces, the so-called
``strange duality'' or ``rank-level duality''. Some evidence
(see \cite{oxbury2} remark 3.4 and \cite{ow}) suggests that this ``strange
duality'' phenomenon should also occur for $\Spin_n$-bundles.

As we have remarked above, the duality was proved in \cite{ram2}
in the case of $\SL_2$ and the proof was later generalized
by the first author to the case of Spin groups. In order to
minimize the  number of papers, we decided to publish our
results together.

It is a pleasure to thank W.M. Oxbury for many useful discussions.

\section{Notation and preliminaries}

\subsection{Prym varieties}
Given a nonzero $\zeta \in JC[2]$, the group of $2$-torsion
points of the Jacobian $JC$, we consider its associated two-sheeted
\'etale cover $\pi_\zeta : \tC_\zeta \ra C$ and the norm map
$\Nm : \pic(\tC_\zeta) \ra \pic(C)$. Mumford \cite{mum1}
introduced the following (isomorphic) subvarieties of $\pic(\tC_\zeta)$,
called Prym varieties
\begin{equation} \label{pryms}
\Nm_\zeta^{-1}(\cO) = P^0_\zeta \cup P'_\zeta,
\qquad  \Nm_\zeta^{-1}(\omega) = P^{ev}_\zeta \cup P^{odd}_\zeta,
\end{equation}
where $P^0_\zeta$ is the connected component containing the origin
$\cO \in JC$  and
$P'_\zeta$ is the other component. The line bundle 
$\omega$ is the canonical line bundle on $C$. For a description of
$\Nm^{-1}_\zeta(\omega)$, see section \ref{exampleparity} Example a). The 
Galois involution $\sigma_\zeta$ of $\tC_\zeta$ acts by
pull-back on the four Prym varieties. We shall also denote the
corresponding involutions by $\sigma_\zeta$. On $P^0_\zeta \cup
P'_\zeta$, we have $\sigma_\zeta(L) = L^{-1}$ and on 
$P^{ev}_\zeta \cup P^{odd}_\zeta$, we have $\sigma_\zeta (\xi) =
\omega \xi^{-1}$.
In this paper we will also use 
$$ P^\zeta_\zeta = \text{one connected component of} \ \Nm^{-1}(\zeta). $$
In order to have a consistent notation, we put 
$P^0_0 = JC$ and $P'_0 = \emptyset$, 
$P^{ev}_0 = \pic^{g-1}(C)$ and $P^{odd}_0 = \emptyset$. On the
finite group $JC[2]$ we have the skew-symmetric Weil pairing which
we denote by $\langle \eta,\zeta \rangle \in \{ \pm 1 \}$. The
group of $2$-torsion points $P^0_\zeta[2]$ of the Prym
$P^0_\zeta$ is isomorphic to the annihilator (w.r.t. the
Weil pairing) of $\zeta \in JC[2]$,i.e., 
\begin{equation} \label{prymtwotorsion}
P^0_\zeta[2] \map{\pi^*_\zeta} \langle \zeta \rangle^\perp / \langle \zeta
\rangle.
\end{equation}
Consider $\eta \in JC[2]$ with $\langle \eta,\zeta \rangle  = 1$. We
also denote by $\eta \in P^0_\zeta[2]$ the corresponding
element under isomorphism \eqref{prymtwotorsion} and
by $T_\eta$ the translation by $\eta$ acting on $P^0_\zeta$ as well
as on the other three Prym varieties \eqref{pryms}.

The Prym variety $P^{ev}_\zeta$ comes equipped with a naturally
defined reduced Riemann theta divisor $\Xi_\zeta$, whose set-theoretical
support equals
\begin{equation} \label{defxi}
\Xi_\zeta := \{ \xi \in P^{ev}_\zeta \ | \ h^0(\tC_\zeta,\xi) >0 \}.
\end{equation}
A translate by a theta-characteristic (resp. a $2$-torsion point)
of $\Xi_\zeta$ gives (non-canonically) a symmetric
theta divisor on $P^0_\zeta$ and $P'_\zeta$ (resp. $P^{odd}_\zeta$)
(see also section 2.2), which we also denote by $\Xi_\zeta$.

Given $\eta \in P^0_\zeta[2]$, we define the line bundle on $P_\zeta$
$$ \LL_\zeta^\eta := \cO_{P_\zeta}(\Xi_\zeta + T_\eta^* \Xi_\zeta). $$
Here $P_\zeta$ stands for any of the four varieties in \eqref{pryms}.
We note that $\Xi_\zeta$ is only canonically defined on $P^{ev}_\zeta$
\eqref{defxi}, but the line bundle 
$\LL_\zeta^\eta$ over any $P_\zeta$ does not 
depend on the choice of $\Xi_\zeta$. Moreover the involutions
$\sigma_\zeta$ of $P_\zeta$ lift canonically to a linear
involution $\sigma_\zeta^*$ of $H^0(P_\zeta,\LL^\eta_\zeta)$. The
subscript $+$ (resp. $-$) will denote the $+$ (resp. $-$) eigenspace
of $\sigma^*_\zeta$.

\subsection{Theta-characteristics} 

In this section we recall some basic results on 
theta-characteristics on Prym varieties. 
Let $\vartheta(C)$ be the set of theta-characteristics on $C$,i.e.,
$\vartheta(C) = \{ \kappa \in \pic^{g-1}(C) \ | \ 
\kappa^{\otimes 2} \map{\sim}
\omega \}$, which comes equipped with a parity map $\epsilon :
\vartheta(C) \ra \{ \pm 1 \}$; $\epsilon(\kappa) = (-1)^{h^0(\kappa)}$.
There is a 1-to-1 correspondence between $\kappa \in \vartheta(C)$
and functions $\tk : JC[2] \ra \{\pm 1\}$ which satisfy $\forall
\eta,\zeta \in JC[2]$
\begin{equation} \label{thetacharquad}
\tk(\eta + \zeta) = \tk(\eta) \tk(\zeta) \langle \eta,\zeta \rangle.
\end{equation}
The correspondence associates to $\kappa \in \vartheta(C)$ the function
$\tk(\eta) = (-1)^{h^0(\kappa \otimes \eta) + h^0(\kappa)}$.
The group $JC[2]$ acts on $\vartheta(C)$ by tensor product and we
have $\epsilon(\eta \kappa) = \epsilon(\kappa) \tk(\eta)$. We define
the set $\vartheta(P^0_\zeta)$ of theta-characteristics on $P^0_\zeta$
to be the set of functions $\tk : P^0_\zeta[2] \cong  
\langle \zeta \rangle^\perp / \langle \zeta \rangle \ra \{ \pm 1\}$
satisfying \eqref{thetacharquad}. Then we have, using
the correspondence between $\kappa$ and $\tk$,
$$ \vartheta(P^0_\zeta) = \{ \kappa \in \vartheta(C) \ | \
\tk(\zeta) = 1 \} / \langle \zeta \rangle $$
Note that we have an equivalence $\kappa \in \vartheta(P^0_\zeta) \iff
\pi^*_\zeta \kappa \in P^{ev}_\zeta$. Given $\kappa \in
\vartheta(P^0_\zeta)$, we consider the symmetric theta divisor
$T_\kappa \Xi_\zeta$ on $P^0_\zeta$. We shall write $\kappa$ instead 
of $\pi^*_\zeta \kappa$. Let $\iota$ be the unique isomorphism
of $\sigma_\zeta^* \cO(T_\kappa \Xi_\zeta)$ with
$\cO(T_\kappa \Xi_\zeta)$ which induces the identity at the origin, then
we have 
\begin{equation} \label{symthetaor}
\iota(\sigma^*_\zeta s_\kappa) = \epsilon(\kappa) s_\kappa
\end{equation}
where $s_\kappa$ is the unique global section of $\cO(T_\kappa \Xi_\zeta)$.
We have similar results for the component $P'_\zeta$. We define
$$ \vartheta(P'_\zeta) = \{ \kappa \in \vartheta(C) \ | \
\tk(\zeta) = -1 \} / \langle \zeta \rangle  = \{ \kappa \in \vartheta(C) \
| \ \pi^*_\zeta \kappa \in P^{odd}_\zeta \}$$
We choose a $2$-torsion point $\alpha \in P'_\zeta[2]$, which allows
us to define uniquely an isomorphism $\iota$ by 
the same condition at $\alpha$. 
Then we still have equality \eqref{symthetaor}.

\subsection{Mumford's parity conservation theorem and Pfaffian line
bundles}

We recall here Mumford's parity conservation theorem \cite{mum2}
which says that the parity of $h^0$ of vector bundles on $C$ is
preserved under deformation, if there exists a non-degenerate 
quadratic form on the family with values in $\omega$.

\subsubsection{Examples} \label{exampleparity}

\noindent
a) $\forall \xi \in P^{ev}_\zeta \cup P^{odd}_\zeta$, the rank $2$
vector bundle $E = \pi_{\zeta*} \xi$ carries a quadratic form
given by the Norm homomorphism $\Nm : \pi_{\zeta*} \xi \ra \omega$.
This explains the notation $P^{ev}_\zeta$ and
$P^{odd}_\zeta$, since the two components of $\Nm^{-1}(\omega)$ may
be distinguished by the property that $h^0(\pi_{\zeta*}\xi)$ is even 
for $\xi \in P^{ev}_\zeta$
and odd for $\xi \in P^{odd}_\zeta$.

\bigskip
\noindent
b) Consider $\eta,\zeta \in JC[2]$ with $\langle \eta,\zeta \rangle =1$.
 $\forall \xi \in P^{ev}_\eta \cup P^{odd}_\eta$ and
$\forall L \in P^0_\zeta \cup P'_\zeta$, the rank $4$ vector bundle
$\pi_{\eta*} \xi \otimes \pi_{\zeta*}L$ has a quadratic form,
which is obtained by multiplying the $\omega$-valued form on 
$\pi_{\eta*} \xi$ and the $\cO$-valued form on $\pi_{\zeta*} L$
(see example a).

\bigskip
\noindent
c) Consider the rank $2n$ vector bundle $\pi_{\zeta*}\xi \otimes
F$, where $\pi_{\zeta*}\xi$ is as in a) and $F$ is an orthogonal 
vector bundle,i.e., equipped with an $\cO$-valued quadratic form.
Multiplying the two forms gives the required $\omega$-valued form
on $\pi_{\zeta*}\xi \otimes F$.

\subsubsection{Pfaffian line bundles} \label{pfaff}

Given a family $\EE$ of vector bundles over $C$ with an
$\omega$-valued quadratic form and parametrized by an integral
variety $S$, one can construct a Pfaffian line bundle $\PF(\EE)$
over $S$ and
a Pfaffian divisor $\div \PF(\EE)$, which are square roots of the
determinant line bundle $\DET(\EE)$ and determinant divisor 
$\div \DET(\EE)$,i.e.,
$$ \DET(\EE) = \PF(\EE)^{\otimes 2} \qquad 2\div \PF(\EE) = \div \DET(\EE) $$
For the construction of $\PF(\EE)$ and $\div\PF(\EE)$ we refer to proposition
7.9 of \cite{ls}.
The support of the Cartier divisor $\div \PF(\EE)$ equals
\begin{equation} \label{supppfaff}
\supp \ \div \PF(\EE)  = \{ s \in S \ | \ h^0(\EE_s) > 0 \}.
\end{equation}
Thus, the Pfaffian divisor associated to the family
$\pi_{\zeta*}\xi$ with $\xi \in P^{ev}_\zeta$ (example a) is the 
Riemann theta divisor $\Xi_{\zeta}$ in \eqref{defxi}.
The Pfaffian divisors associated to the families of example b) (resp. 
example c) will appear in the proof of proposition \ref{pwd} (resp.
in the construction of the Prym-Spin duality, see
section 3).

\section{Prym varieties and the moduli of $\Spin_n$-bundles}

In this section we recall the definition of the linear maps of 
theorem \ref{duality} (see \cite{oxbury1},\cite{oxbury2}).
Let $\mathrm{SC}_n$ be the special Clifford group of a 
non-degenerate quadratic form of dimension $n$ and $\MM(\mathrm{SC}_n)$
be the moduli space of semistable principal $\mathrm{SC}_n$-bundles. The
spinor norm induces a surjective morphism $\Nm : \MM(\mathrm{SC}_n)
\ra \pic(C)$ and we define 
$$ \MM^+(\Spin_n) = \MM(\Spin_n) = \Nm^{-1}(\cO),
\qquad
\MM^-(\Spin_n) = \Nm^{-1}(\cO(p)), $$
where $p$ is a fixed point on $C$. In section 6 we shall give  
examples of the moduli spaces $\MM^+(\Spin_n)$ for low values of $n$.
We denote by
$V \cong \cc^n$ the standard orthogonal representation of
$\mathrm{SC}_n$ and by $S$ (if $n$ is odd) and $S^\pm$ (if $n$ is even) the
spinor and half-spinor representations. Given a Clifford or Spin-bundle $A$,
we shall denote by $A(V),A(S),A(S^\pm)$ the induced vector bundles.

There is a natural action of $JC[2]$ on the moduli spaces 
$\MM^\pm(\Spin_n)$ and we recall that there are natural Galois
covers 
\begin{equation} \label{galcovso}
\MM^\pm(\Spin_n) \map{JC[2]} \MM^\pm(\SO_n), \qquad A \ra A(V).
\end{equation}
The Galois action of $JC[2]$ on the moduli $\MM^\pm(\Spin_n)$ is given 
by $A \mapsto \alpha . A$ for $\alpha \in JC[2]$ with 
\begin{equation} \label{galaction}
(\alpha . A)(V) = A(V) \qquad (\alpha . A)(S^\pm) = A(S^\pm) \otimes \alpha
\end{equation}
The two components $\MM^\pm(\SO_n)$ of the moduli of semistable 
orthogonal bundles are distinguished by the second Stiefel-Whitney
class of the bundles. For some general facts on orthogonal bundles,
see \cite{ramanan}.

To construct the linear maps  of theorem \ref{duality}, it will
be enough to exhibit effective Cartier divisors in the right linear
systems
$$\DD_{2m+1,\eta}^\pm \subset P^{ev}_\eta \times \MM^\pm(\Spin_{2m+1}),
\qquad 
\DD_{2m,\eta}^\pm \subset (P^{ev}_\eta \cup P^{odd}_\eta) 
\times \MM^\pm(\Spin_{2m}).$$
We already observed in \ref{exampleparity} example c) that the
family $\EE : = \{ \pi_{\eta*} \xi \otimes A(V) \}$ carries an
$\omega$-valued form, so we define $\DD_{n,\eta} := \div \PF(\EE)$. 
By \eqref{supppfaff} we have (e.g. for $n$ odd)
\begin{equation} \label{defDeta}
\supp \  \DD^\pm_{2m+1,\eta} = \{ (\xi,A) \in P^{ev}_\eta \times
\MM^\pm(\Spin_{2m+1}) \ | \ h^0(C, \pi_{\eta*} \xi \otimes 
A(V)) >0 \}.
\end{equation}
We shall use the abreviated notation $\DD_\eta$, if there is
no ambiguity. We notice that the divisor $\DD_\eta$ does not
descend to a Cartier divisor on $P_\eta \times \MM^\pm(\SO_n)$.
We denote by $s^\pm_\eta$ their associated linear maps. 
The details of this construction are worked out in \cite{oxbury2}
section 6.

Given a pair of semistable orthogonal vector bundles $(E,E')$ in
$\MM^\pm(\SO_n) \times \MM^+(\SO_{n'})$ we can associate its
orthogonal direct sum $E \oplus E' \in \MM^\pm(\SO_{n+n'})$, which gives rise
to a morphism $\MM^\pm(\SO_n) \times \MM^+(\SO_{n'}) \ra
 \MM^\pm(\SO_{n+n'})$. In the next lemma we will show that this
morphism lifts to the moduli spaces of $\Spin$-bundles.

\begin{lem} \label{addspinbundles}
For any $n,m \geq 1$, there are natural morphisms
\begin{equation} \label{evenaddspin}
1. \ \ \ \ \MM^\pm(\Spin_{2n}) \times \MM^+(\Spin_{2m}) \map{\iota}
\MM^\pm(\Spin_{2n+2m}) 
\end{equation}
\begin{equation} \label{oddaddspin}
2. \ \ \ \ \MM^\pm(\Spin_{2n+1}) \times \MM^+(\Spin_{2m}) 
\map{\iota} \MM^\pm(\Spin_{2n+2m+1}) 
\end{equation}
which are lifts via \eqref{galcovso} of the direct sum morphisms
for orthogonal vector bundles. Given a pair $(A,A')$ of $\Spin$-bundles,
then its sum $A+A' := \iota(A,A')$ has the following properties
\begin{enumerate}
\item the two associated half-spinor vector bundles are 
$$(A+A')(S^+) = A(S^+) \otimes A'(S^+) \oplus A(S^-) \otimes A'(S^-)$$
$$(A+A')(S^-) = A(S^+) \otimes A'(S^-) \oplus A(S^-) \otimes A'(S^+)$$
\item the associated spinor vector bundle is 
$$(A+A')(S) = A(S) \otimes (A'(S^+) \oplus A'(S^-))$$
\end{enumerate}
Moreover the pull-back by $\iota$ of the determinant line bundle
decomposes ($l=2n$ or $2n+1$)
\begin{equation} \label{decthetalb}
\iota^* \Theta(\cc^{l+2m}) = \Theta(\cc^l) \boxtimes \Theta(\cc^{2m}) 
\end{equation}
\end{lem}

\begin{proof} 
The natural homomorphism of algebraic groups 
$\mathrm{SC}_l \times \Spin_{2m} \ra \mathrm{SC}_{l+2m}$ induces a
morphism at the level of Clifford (and Spin)-bundles. 
By \cite{ramanan} prop.4.2
and \cite{oxbury1} lemma 1.2 semistability is preserved under
this operation. Thus, taking into account that the spinor
norm is preserved under the above group homomorphism and that
$\MM^\pm(\Spin_n)$ is a coarse moduli space, we get the above
claimed morphisms. The expressions of the associated (half-)spinor
bundles are easily deduced from the definitions of the (half-)spinor
representations of the Spin groups.
\end{proof}

\section{Multiplication of theta functions}
In this section we recall some facts on theta functions over an 
abelian variety. Let $(A,\Theta)$ be a principally polarized
abelian variety, where $\Theta$ is
a symmetric divisor representing the polarization. The subscript
$\pm$ denotes the $\pm$eigenspaces of $H^0(A,m\Theta)$ under the
canonical involution of $A$. Recall that $H^0(A,2\Theta)_+ =
H^0(A,2\Theta)$. We will need
the following facts on multiplication maps 

\begin{prop} \label{multkk}
If $l \geq2$ and $m \geq 3$, then the multiplication maps
\begin{enumerate} 
\item 
$H^0(A,l\Theta) \otimes H^0(A,m\Theta) \lra H^0(A,(l+m)\Theta) $
\item
$H^0(A,2l\Theta)_+ \otimes H^0(A,m\Theta) \lra H^0(A,(2l+m)\Theta)$
\end{enumerate}
are surjective.
\end{prop}

For a proof of 1, we refer e.g. to \cite{kempf}. For a proof of 2, see
prop. 1.4.4 \cite{khaled}.

\bigskip
The symmetric theta divisor $\Theta$ allows us to identify $A$ with its
dual variety $\hat{A}$. Let $m$ be the isogeny
\begin{eqnarray} \label{isogeny}
 m \ : \ A \times A & \lra & A \times A \\
  (x,y) & \lms & (x+y,x-y) \nonumber
\end{eqnarray}
then it is well-known that for any $\alpha \in A[2]= \hat{A}[2]$, we have 
$$ m^* \cO_A(2\Theta \otimes \alpha) \boxtimes
\cO_A(2\Theta \otimes \alpha) = \cO_A(4\Theta) \boxtimes
\cO_A(4\Theta) $$
We take global sections and take the sum over all $\alpha \in
A[2]$ to get the direct sum decomposition
\begin{equation} \label{dsdm}
m^*: \sum_{\alpha \in A[2]}  H^0(A,2\Theta \otimes \alpha)
\otimes H^0(A,2\Theta \otimes \alpha)
\map{\sim}  H^0(A,4\Theta) \otimes  H^0(A,4\Theta)
\end{equation}
Since $m$ is equivariant for the canonical involution of $A \times
A$, we obtain the following decomposition into
$\pm$eigenspaces
\begin{equation} \label{iso4+}
\sum_{\alpha \in A[2]} H^0_+ \otimes H^0_+ \oplus  H^0_- \otimes H^0_- 
= \left[  H^0(A,4\Theta) \otimes  H^0(A,4\Theta) \right]_+ 
\end{equation}
\begin{equation} \label{iso4-}
\sum_{\alpha \in A[2]} H^0_+ \otimes H^0_- \oplus H^0_- \otimes H^0_+ =
\left[  H^0(A,4\Theta) \otimes  H^0(A,4\Theta) \right]_-
\end{equation}
In particular, if we restrict $m$ to the diagonal $A \inj A \times A$
we get the direct sum decomposition (in this case we get the 
duplication map $m_{|A} = [2] :\ A \ra A$ ; $x \mapsto 2x$)
\begin{equation} \label{dupl}
[2]^* : \sum_{\alpha \in A[2]}  
H^0(A,2\Theta \otimes \alpha) \map{\sim}
H^0(A,8\Theta)
\end{equation}
In the next proposition we will need a more general version, $\forall
n \geq 1$
\begin{equation} \label{dupl2}
[2]^* : \sum_{\alpha \in A[2]}  
H^0(A,n\Theta \otimes \alpha) \map{\sim}
H^0(A,4n\Theta)
\end{equation}

Via the Weil pairing we can identify $A[2]$ with the group 
of characters $\Hom(A[2],\cc^*)$; $\alpha \mapsto \langle \alpha, \cdot
\rangle$. The natural action of the group $A[2]$ on $A$ and
$A \times A$ lifts canonically to a linear action on the RHS spaces
of \eqref{dsdm}, \eqref{dupl} and \eqref{dupl2}. Then the
direct sum decompositions of \eqref{dsdm}, \eqref{dupl} and \eqref{dupl2}
are precisely the character space decompositions of this 
group action.

\begin{prop} \label{multmap}
The following multiplication maps are surjective
\begin{enumerate}
\item
$\sum_{\alpha \in A[2]} H^0(A,2\Theta \otimes \alpha)_\pm
\otimes H^0(A,\Theta \otimes \alpha) \lra H^0(A,3\Theta)_\pm$
\item
$\sum_{\alpha \in A[2]} H^0(A,2\Theta \otimes \alpha)_+
\otimes H^0(A,2\Theta \otimes \alpha)_\pm \lra H^0(A,4\Theta)_\pm$
\end{enumerate}
\end{prop}

\begin{proof}
First let us prove 1. By proposition \ref{multkk}.1 ($l=4,m=8$), the
multiplication map
\begin{equation} \label{multmap48}
H^0(A,4\Theta) \otimes H^0(A,8\Theta) \lra H^0(A,12\Theta)
\end{equation}
is surjective. Consider the character space decomposition
\eqref{dupl2} of the two spaces. Since the above tensor product
is compatible with the linear action of the group $A[2]$, we get 
a surjective map between the two character spaces of \eqref{multmap48}
corresponding to the zero character,i.e.,
$$\sum_{\alpha \in A[2]} H^0(A,2\Theta \otimes \alpha)
\otimes H^0(A,\Theta \otimes \alpha) \lra H^0(A,3\Theta)$$
Considering $\pm$eigenspaces proves assertion 1. By proposition 
\ref{multkk}.2
($l=4,m=8$), the multiplication map 
$$H^0(A,8\Theta)_+ \otimes H^0(A,8\Theta) \lra 
H^0(A,16\Theta)$$
is surjective. Considering the zero character space and $\pm$eigenspaces, we
will prove 2. 
\end{proof}

Finally, let us denote by $\Heis$ the Heisenberg group associated
to the line bundle $\cO(2\Theta)$ (see \cite{mum3} or 
\cite{beauville} page 280),i.e. a central extension of $A[2]$ by
$\cc^*$
\begin{equation} \label{heisenberg}
0 \lra \cc^* \lra \Heis \lra A[2] \lra 0 
\end{equation}
We recall that $\Heis$ acts linearly on $H^0(A,2\Theta)$, its
unique (up to conjugation) representation of level $1$.

\section{The Prym-Wirtinger duality}

Consider $\eta,\zeta \in JC[2]$ such that $\langle \eta,\zeta \rangle = 1$.
We recall the definition of the line bundle $\LL_\zeta^\eta = \cO_{P_\zeta}
(\Xi_\zeta + T^*_\eta \Xi_\zeta) = \cO_{P_\zeta}(2\Xi_\zeta \otimes
\eta)$. In the next proposition we show that exchanging the roles
of $\eta$ and $\zeta$ will establish a duality at the level of global
sections of $\LL^\eta_\zeta$. We may view this duality
as an analogue for Prym varieties of the well-known Wirtinger
duality (put $\eta = \zeta = 0$) for Jacobians (see \cite{mum1}
page 335). This Prym-Wirtinger duality is at the heart of the
proof of theorem 1.1 (section 7), as the Prym-Spin pairing
``restricts'' to  the Prym-Wirtinger duality for suitably chosen (products 
of) Prym varieties in $\MM(\Spin_n)$.

\begin{prop} \label{pwd}
We have the following canonical isomorphisms for any $\eta,\zeta \in
JC[2]$ such that $\langle \eta,\zeta \rangle = 1$
\begin{equation} \label{pwd+}
 H^0(P^{ev}_\eta, \LL_\eta^\zeta)^\vee_+ 
  \cong   H^0(P^0_\zeta,\LL_\zeta^\eta)_+  
\qquad
 H^0(P^{odd}_\eta, \LL_\eta^\zeta)^\vee_+ 
  \cong   H^0(P^0_\zeta, \LL_\zeta^\eta)_-  
\end{equation}
\begin{equation} \label{pwd-}
 H^0(P^{ev}_\eta, \LL_\eta^\zeta)^\vee_- 
  \cong   H^0(P'_\zeta, \LL_\zeta^\eta)_+  
\qquad
 H^0(P^{odd}_\eta, \LL_\eta^\zeta)^\vee_- 
  \cong   H^0(P'_\zeta, \LL_\zeta^\eta)_-
\end{equation}
\end{prop}

\begin{proof}
We will show that the duality between the above vector spaces
is given by a reduced Cartier divisor, whose set-theoretical
support equals
\begin{equation} \label{divpwd}
\Delta^{Wirt}_{\eta,\zeta} := \{ (\xi,L) \in
(P^{ev}_\eta  \cup P^{odd}_\eta) \times (P^0_\zeta \cup
P'_\zeta) \ | \  h^0(C,\pi_{\eta*}\xi \otimes \pi_{\zeta*}L) >0 \}.
\end{equation}
Indeed, as shown in example b) of \ref{exampleparity}, the family
$\EE$ of rank $4$ vector bundles $\pi_{\eta*}\xi \otimes 
\pi_{\zeta*} L$ parametrized by $(P^{ev}_\eta  \cup P^{odd}_\eta) \times
(P^0_\zeta \cup P'_\zeta)$ is equipped with an $\omega$-valued
quadratic form. Hence, by \ref{pfaff}, we can consider its 
associated Pfaffian divisor  $\Delta^{Wirt}_{\eta,\zeta} : = 
\div \PF (\EE)$. By \eqref{supppfaff} its support equals the set in
\eqref{divpwd} and an easy computation shows that 
$\cO(\Delta^{Wirt}_{\eta,\zeta}) = \LL^\zeta_\eta \boxtimes 
\LL^\eta_\zeta$. Hence we obtain a pairing
\begin{equation} \label{pairingpw}
H^0(P^{ev}_\eta  \cup P^{odd}_\eta, \LL^\zeta_\eta)^\vee \map{\psi}
H^0(P^0_\zeta \cup P'_\zeta, \LL^\eta_\zeta).
\end{equation}
First we will show that $\psi$ is an isomorphism. We
verify that both sides of \eqref{pairingpw} are $\Heis$-modules
\eqref{heisenberg}
of level $1$ and that $\psi$ is $\Heis$-equivariant. Hence, since
$\psi$ is nonzero, it is an isomorphism.

To finish the proof we have to analyze how $\psi$ acts on the
$\pm$eigenspaces of the linear involutions $\sigma^*_\eta$ (resp.
$\sigma^*_\zeta$). Let us restrict attention to the duality on the
component $P^{ev}_\eta \times P^0_\zeta$. Consider the rational map,
induced by the divisor $\Delta^{Wirt}_{\eta,\zeta}$
\begin{eqnarray*}
\Delta : P^{ev}_\eta & \lra & \pp H^0(P^0_\zeta, \LL^\eta_\zeta) \\
 \xi & \lms & \Delta(\xi) := {\Delta^{Wirt}_{\eta,\zeta}}_{|\{\xi\} \times
P^0_\zeta}
\end{eqnarray*}
We observe that, $\forall \xi \in P^{ev}_\eta$, the divisor $\Delta(\xi)$
is invariant under the (projective) involution $\sigma^*_\zeta$. Consider
$\kappa \in \vartheta(C)$ such that $\kappa \in \vartheta(P^0_\zeta) \cap
\vartheta(P^0_\eta)$,i.e., $\tk(\eta) = \tk(\zeta) = 1$. Then it follows
from the definition of $\Delta^{Wirt}_{\eta,\zeta}$ that
$\Delta(\pi^*_\eta \kappa) = T^*_\kappa \Xi_\zeta +
T^*_{\eta \kappa} \Xi_\zeta$. By \eqref{symthetaor} and with the
notation as above, we have 
$$ \iota \sigma^*_\zeta (s_\kappa \cdot s_{\eta \kappa}) = 
\epsilon(\kappa) \epsilon(\eta \kappa) s_\kappa \cdot s_{\eta \kappa} =
\tk(\eta) s_\kappa \cdot s_{\eta \kappa} = s_\kappa \cdot s_{\eta \kappa} $$
So we get $\Delta(\pi^*_\eta \kappa) \in \pp H^0(P^0_\zeta, \LL^\eta_\zeta)_+$.
Hence $\im \Delta \subset \pp H^0(P^0_\zeta, \LL^\eta_\zeta)_+$. On the
other hand, we see that the divisor $\Delta(L) :=
{\Delta^{Wirt}_{\eta,\zeta}}_{| P^{ev}_\eta \times \{L\}}$
for $L \in P^0_\zeta$ is symmetric and, take e.g. $L=\cO$, $\Delta(L) \in
\pp H^0(P^{ev}_\eta, \LL^\zeta_\eta)_+$. Hence the pairing 
\eqref{pairingpw} splits as follows
$$H^0(P^{ev}_\eta, \LL^\zeta_\eta)_+^\vee \map{\sim}   
H^0(P^0_\zeta, \LL^\eta_\zeta)_+ $$
By the same  reasoning as above, one gets the other three isomorphisms as
stated in the proposition. We leave the details to the reader.
\end{proof}

\section{The Schottky variety $SCH_n$}

In this section we are going to define maps from products of
Prym varieties into $\MM^{\pm}(\Spin_n)$. The disjoint union of these
products will be called the Schottky variety, denoted by
$SCH_n$. Let us start with the first cases $n \leq 6$, which
will serve as a pattern to construct by induction the Schottky  
variety for any integer $n$.

\subsection{$n=2$} \label{spin2}

This case has been extensively discussed in \cite{oxbury2}
section 4.2.
Unlike in \cite{oxbury2}, we will define $\MM^+(\Spin_2)$ to be the
degree $0$ component of $\MM(\cc^*) = \pic(C)$,i.e., $\MM^+(\Spin_2) =
JC$. We also put $\MM^-(\Spin_2) = \emptyset$. 
We define the Schottky variety by $SCH_2 :=JC$ and the map
$\phi_2$ is the identity $\phi_2 : JC = \MM^+(\Spin_2)$. 
Given $A \in \MM^+(\Spin_2)$, we denote the 
line bundle $A(S^+) \in JC$ by $N$. Then we have $A(S^-)= N^{-1}$
and $A(V) = N^2 \oplus N^{-2}$.

\subsection{$n=3$} \label{spin3}

In this case (see \cite{oxbury2} section 4.3) we have the isomorphisms
$\MM^+(\Spin_3) = \mod0$ and $\MM^-(\Spin_3) = \modp$,
where we used the exceptional isomorphism
$\Spin_3 \cong \SL_2$ at the level of algebraic groups. 
Here $\mathcal{SU}_C(2,L)$ denotes the moduli space of semistable
rank $2$ vector bundles with fixed determinant equal to $L \in \pic(C)$.
To  a semistable $\Spin_3$-bundle $A \in \MM^+(\Spin_3)$ 
we associate two semistable vector bundles $A(S)$ (resp. $A(V)$)
induced by the spinor (resp. orthogonal) representation.  
In particular the isomorphism $\MM^+(\Spin_3) \map{\sim} \mod0$ is given
by the map $A \mapsto A(S)$.
Then the orthogonal 
rank $3$ bundle is given by
$$ A(V) = \End_0(A(S)). $$
We choose a point $\tp \in \tC_\zeta$ lying over $p$.
Taking direct image gives the following maps
$$
\begin{array}{rcl}
\phi^+_{3,\zeta} :  P^\zeta_\zeta & \lra & \MM^+(\Spin_3) \\
 L & \lms & \pi_{\zeta*}L = A(S)
\end{array}
\qquad
\begin{array}{rcl}
\phi^-_{3,\zeta} :  P^\zeta_\zeta & \lra & \MM^-(\Spin_3) \\
 L & \lms & \pi_{\zeta*} L(\tp) = A(S).
\end{array}
$$
If $\zeta =0$, we define $\phi^+_{3,0}(L) = L \oplus L^{-1}$ and
$\phi^-_{3,0}(L)$ is the unique stable rank $2$ bundle, which 
fits into the exact sequence (see \cite{beauville} section 3)
$$ 0 \lra L \oplus L^{-1} \lra \phi^-_{3,0}(L) \lra \cc_p \lra 0. $$
\begin{lem} \label{lemdec3}
Given $M \in \pic(\tC_\zeta)$. Then the orthogonal bundle
$\End_0(\pi_{\zeta*}M)$ is an orthogonal direct sum
$$ \End_0(\pi_{\zeta*}M) = \zeta \oplus \pi_{\zeta*} \Big( M^2 \otimes
\pi^*_\zeta(\Nm M)^{-1} \Big). $$
\end{lem}

\begin{proof}
Since $\pi^*_{\zeta} \zeta = \cO_{\tC}$, we have an isomorphism
$\zeta \otimes \pi_{\zeta*}M \map{\sim} \pi_{\zeta*}M$, which
gives rise to a homomorphism $\zeta \ra \End_0(\pi_{\zeta*}M)$.
Since both $\zeta$ and $\End_0(\pi_{\zeta*}M)$ have degree $0$ and
since $\End_0(\pi_{\zeta*}M)$ is poly-stable, it follows that
$\zeta$ is a subbundle of $\End_0(\pi_{\zeta*}M)$. It remains to
determine its supplement bundle. There is a natural homomorphism
of $\sym^2(\pi_{\zeta*}M)$ into $\pi_{\zeta*}M^2$. Using the 
isomorphism $\End_0(\pi_{\zeta*}M) \cong \sym^2(\pi_{\zeta*}M) \otimes
(\Nm M)^{-1} \otimes \zeta$, we get a  homomorphism of 
$\End_0(\pi_{\zeta*} M)$ into $\pi_{\zeta*}M^2 \otimes (\Nm M)^{-1}$.
A pointwise check shows that this is a surjective homomorphism
and that it is supplementary to $\zeta$.
\end{proof}

\noindent
In particular, we have for $L \in P^\zeta_\zeta$
\begin{equation} \label{decend}
\End_0(\pi_{\zeta*}L) = \zeta \oplus \pi_{\zeta*} L^2 
\qquad
\End_0(\pi_{\zeta*}L(\tp)) = \zeta \oplus \pi_{\zeta*} (
L^2(\tp - \sigma \tp))
\end{equation}

\noindent
We define the Schottky variety and the morphism $\phi^\pm_3$
$$SCH_3 : = \coprod_{\zeta \in JC[2]} P^\zeta_\zeta, \qquad
\phi^\pm_3 := \coprod_{\zeta} \phi^\pm_{3,\zeta} : \ SCH_3 \lra
 \MM^\pm(\Spin_3)$$

\subsection{$n=4$} \label{spin4}

In this case (see \cite{oxbury2} section 4.4) we have the
isomorphisms $\MM^+(\Spin_4) = \mod0 \times \mod0$ and
$\MM^-(\Spin_4) = \modp \times \modp$, where we
used the exceptional isomorphism $\Spin_4 = \SL_2 \times \SL_2$.
The previous isomorphism is given by sending a $\Spin_4$-bundle
$A$ to the pair of rank $2$ bundles $(A(S^+),A(S^-))$. Consider
the maps
\begin{eqnarray*}
 \phi^\pm_{4,\zeta} : P^\zeta_\zeta \times P^\zeta_\zeta &  \lra & 
\MM^\pm(\Spin_4) \\
(L,M) & \lms & \phi^\pm_{4,\zeta}(L,M) := (\phi_{3,\zeta}^\pm(L),
\phi_{3,\zeta}^\pm(M))
\end{eqnarray*}
where we use the previous isomorphisms.
\begin{lem} \label{lemdec4}
Given a pair $(L,M) \in P^\zeta_\zeta \times P^\zeta_\zeta$,
we consider their associated $\Spin_4$-bundles $A^\pm :=
\phi^\pm_4(L,M)$. Then we have
$$ A^+(V) = \pi_{\zeta*}(LM) \oplus \pi_{\zeta*}
(LM^{-1}), \qquad
A^-(V) = \pi_{\zeta*}(LM(\tp -\sigma \tp)) \oplus \pi_{\zeta*}
(LM^{-1}). 
$$
\end{lem}

\begin{proof}
We know that the orthogonal bundle $A^\pm(V)$ is the tensor
product $A^\pm(S^+) \otimes A^\pm(S^-)^\vee$. Let us do the
computations for $A^+(V)$. 
$$
A^+(V) = A^+(S^+) \otimes A^+(S^-) = \pi_{\zeta*}L \otimes 
\pi_{\zeta*}M = \pi_{\zeta*}(L \otimes \pi_\zeta^* \pi_{\zeta*}M)
$$
where we used the fact that $A^+(S^-)$ is self-dual and the projection
formula for the map $\pi_{\zeta}$. Now by example a) 
\ref{exampleparity}, $\pi_\zeta^* \pi_{\zeta*}M$ is a semistable
orthogonal bundle over $\tC_\zeta$, hence this bundle splits
$\pi_\zeta^* \pi_{\zeta*}M = M \oplus M^{-1}$. The computations
for $A^-(V)$ are similar. 
\end{proof}

\noindent
We define the Schottky variety and the morphism $\phi^\pm_4$ by
$$\coprod \phi^\pm_{4,\zeta}: SCH_4 : = \coprod_{\zeta \in JC[2]} 
P^\zeta_\zeta \times P^\zeta_\zeta \lra \MM^\pm(\Spin_4).$$

\subsection{$n=5$}

We define $\phi^\pm_{5,\zeta}$ to be the composite map (see sections
\ref{spin2}, \ref{spin3} and lemma \ref{addspinbundles})
$$P^\zeta_\zeta \times JC \map{\phi^\pm_{3,\zeta} \times \phi_2}
\MM^\pm(\Spin_3) \times \MM^+(\Spin_2) \map{\iota} \MM^\pm(\Spin_5) $$
For example, given $(L,N) \in P^\zeta_\zeta \times JC$, the $\Spin_5$-bundle
$A:= \phi^+_{5,\zeta}(L,N)$ satisfies 
\begin{equation} \label{lemdec5}
A(S) = (N\oplus N^{-1}) \otimes \pi_{\zeta*}L \qquad
A(V) = \zeta \oplus N^2 \oplus N^{-2} \oplus \pi_{\zeta*}L^2 
\end{equation}
We remark that in this case we have an isomorphism $\MM^+(\Spin_5) =
\MM^+(\mathrm{Sp}_4), A \mapsto A(S)$. The symplectic form 
on the bundle  $A(S)$ in \eqref{lemdec5} is the obvious one. We define
the Schottky variety by
$$ \phi_5^\pm : SCH_5 : = \coprod_{\zeta \in JC[2]} P^\zeta_\zeta \times
JC \lra \MM^\pm(\Spin_5)$$

\subsection{$n=6$}

Consider the composite map (see sections \ref{spin2}, \ref{spin4} and
lemma \ref{addspinbundles})
$$ \phi^\pm_{6,\zeta} : P^\zeta_\zeta \times P^\zeta_\zeta \times JC
 \map{\phi^\pm_{4,\zeta} \times \phi_2} \MM^\pm(\Spin_4) \times \MM^+(\Spin_2) 
\map{\iota} \MM^\pm(\Spin_6) $$
which associates to the triple $(L,M,N) \in
 P^\zeta_\zeta \times P^\zeta_\zeta \times JC$ the $\Spin_6$-bundle
$A$, which verifies
$$A(V) = N^2 \oplus N^{-2} \oplus \pi_{\zeta*}LM \oplus
\pi_{\zeta*}LM^{-1}$$ 
$$ A(S^+) = (N \otimes \pi_{\zeta*}M) \oplus 
(N^{-1} \otimes \pi_{\zeta*}L)
\qquad 
A(S^-) = (N \otimes \pi_{\zeta*}L) \oplus (N^{-1} \otimes \pi_{\zeta*}M) $$
We define the Schottky variety by
$$\phi_6^\pm: SCH_6 : = \coprod_{\zeta \in JC[2]} P^\zeta_\zeta \times
P^\zeta_\zeta \times JC \lra \MM^\pm(\Spin_6)$$

\subsection{the general case}

Now we will define by induction the Schottky variety $SCH_n$ for any
integer. We put for $n \geq 3$
$$ SCH_{n+4,\zeta} := P^\zeta_\zeta \times P^\zeta_\zeta
\times SCH_n \qquad  SCH_{n+4} := \coprod_{\zeta \in JC[2]} 
SCH_{n+4,\zeta} $$
and define on each component $SCH_{n+4,\zeta}$ the morphism 
$\phi^\pm_{n+4,\zeta}$ as the composite map (see lemma \ref{addspinbundles})
$$  SCH_{n+4,\zeta}
\map{\phi^+_{4,\zeta} \times \phi^\pm_n} \MM^+(\Spin_4) 
\times \MM^\pm(\Spin_n)
\map{\iota} \MM^\pm(\Spin_{n+4}) $$
Note that for $n \geq 5$ the image of $SCH_n$ is contained in the
semistable boundary of $\MM^\pm(\Spin_n)$.

\section{Proof of theorem 1.1}
As in the previous section, we first will prove the
cases $n \leq 6$ separately and then proceed by induction to prove the 
duality for any integer $n$.

\subsection{$n=2$}
For the details see \cite{oxbury2} section 4.2. In this case
$\MM^+(\Spin_2) = JC$ and $\Theta(\cc^2) = \cO_{JC}(8\Theta)$.
The morphism \eqref{prymspineven} decomposes as follows:
by the Prym-Wirtinger duality \eqref{pwd+} with $\zeta =0$, we
have 
$$ \sum_{\eta \in JC[2]} H^0(P^{ev}_\eta, 2\Xi_\eta)^\vee \oplus
H^0(P^{odd}_\eta, 2\Xi_\eta)^\vee  = 
\sum_{\eta \in JC[2]} H^0(JC,2\Theta\otimes \eta)_+ \oplus
  H^0(JC,2\Theta \otimes \eta)_-, $$
which by \eqref{dupl} equals $H^0(JC,8\Theta)$.

\subsection{$n=3$}

First we shall consider the duality on $\MM^+(\Spin_3)$.
By the first equality in \eqref{lemdec3}, we have 
$ \forall L \in P^\zeta_\zeta, \ \forall \xi \in
P^{ev}_\eta$ 
\begin{equation} \label{dech0}
 h^0(C,\End_0(\pi_{\zeta*}L) \otimes
\pi_{\eta*}\xi) = h^0(\tC_\eta, \xi\otimes \pi_\eta^* \zeta)
+ h^0(C,\pi_{\zeta*}L^2 \otimes \pi_{\eta*}\xi)
\end{equation}
Hence we claim that the pull-back of the divisor $\DD^+_\eta$
\eqref{defDeta} by the morphism
$$\id \times \phi^+_{3,\zeta}: \ P^{ev}_\eta \times P^\zeta_\zeta 
\lra P^{ev}_\eta \times \MM^+(\Spin_3) $$
splits into two divisors
\begin{equation} \label{decdiveta3}
(\id \times \phi^+_{3,\zeta})^*(\DD^+_\eta) = p_1^*(T^*_{\zeta}
\Xi_\eta) + (\id \times [2])^* (\Delta^{Wirt}_{\eta,\zeta})
\end{equation}
where $p_1$ is the projection onto the first factor $P^{ev}_\eta$, $[2]: 
P^\zeta_\zeta \ra P^0_\zeta$ is the duplication map $L \mapsto L^2$,
and $\Delta^{Wirt}_{\eta,\zeta}$ is the divisor defined in
\eqref{divpwd}. Indeed, by \eqref{dech0} the decompositon \eqref{decdiveta3}
follows set-theoretically and since $(\id \times \phi^+_{3,\zeta})^*
(\DD^+_\eta) \in |3\Xi_\eta \boxtimes 8 \Xi_\zeta|$, the equality
\eqref{decdiveta3} also holds scheme-theoretically. The next lemma
is an immediate consequence of this decomposition.

\begin{lem} \label{lemfact3}
For any $\eta,\zeta$ satisfying $\langle \eta,\zeta \rangle = 1$,
the linear map $s^+_\eta$ composed with the pull-back induced
by $\phi^+_{3,\zeta}$
$$H^0(P^{ev}_\eta,3\Xi_\eta)^\vee_+ \map{s^+_\eta}
H^0(\MM^+(\Spin_3), \Theta(\cc^3)) \map{\phi^{+*}_{3,\zeta}}
H^0(P^\zeta_\zeta,8\Xi_\zeta)_+$$
factorizes as follows
$$H^0(P^{ev}_\eta,3\Xi_\eta)^\vee_+ \map{+T^*_\zeta \Xi_\eta}
H^0(P^{ev}_\eta,\LL_\eta^\zeta)^\vee_+
\map{\eqref{pwd+}} H^0(P^0_\zeta,\LL_\zeta^\eta)_+ \map{[2]^*} 
H^0(P^\zeta_\zeta, 8\Xi_\zeta)_+$$
where the first map is the dual of the multiplication map $D \mapsto
D + T^*_\zeta \Xi_\eta$.
\end{lem}

We are now in a position to prove \eqref{prymspinodd} for $m=1$.
The main idea of the proof is to show that the composite map
\begin{equation} \label{dualsch3}
\sum_\eta H^0(P^{ev}_\eta,3\Xi_\eta)^\vee_+
\map{\sum s^+_\eta} H^0(\MM^+(\Spin_3), \Theta(\cc^3))
 \map{\phi^{+*}_3} \sum_\zeta H^0(P^\zeta_\zeta,8\Xi_\zeta)_+
\end{equation}
is injective, which immediately implies that $\sum s^+_\eta$ 
is an isomorphism, since the first two
spaces have the same dimension (\cite{oxbury2} theorem 3.1).

By lemma \ref{lemfact3}, the linear map \eqref{dualsch3} factorizes
as follows
\begin{eqnarray*}
\sum_{\eta \in JC[2]} H^0(P^{ev}_\eta,3\Xi_\eta)^\vee_+ & \map{(1)} &
\sum_{\eta \in JC[2]} \sum_{\zeta \in P^0_\eta[2]}
H^0(P^{ev}_\eta,\LL_\eta^\zeta)^\vee_+ \\
 & \map{(2)} & \sum_{\zeta \in JC[2]} \sum_{\eta \in P^0_\zeta[2]}
  H^0(P^0_\zeta, \LL_\zeta^\eta)_+ \\
 & \map{(3)} & \sum_{\zeta \in JC[2]} 
 H^0(P^\zeta_\zeta,8\Xi_\zeta)_+ 
\end{eqnarray*}
The arrows are as follows: (1) is the dual of the multiplication map
in proposition \ref{multmap}.1, which is injective; (2) is the 
Prym-Wirtinger duality \eqref{pwd+}, which is an isomorphism and (3)  is the
isomorphism induced by the duplication map \eqref{dupl}. Finally
we observe that we can invert the indices of summation in (2), since
both sets of indices are in 1-to-1 correspondence with 
the set of isotropic (w.r.t
the Weil form) Klein subgroups of $JC[2]$. Since the three maps (1),(2),(3)
are injective, their composite map \eqref{dualsch3} is also
injective and we are done.

\bigskip
Let us briefly indicate how to adapt the previous proof to the
moduli $\MM^-(\Spin_3)$. First, using the second equality in 
\eqref{lemdec3}, we easily see that the analogue of 
\eqref{decdiveta3} is
$$ (\id \times \phi^-_{3,\zeta})^*(\DD^-_\eta) = p_1^*(T^*_{\zeta}
\Xi_\eta) + (\id \times (T_{\tp}\circ [2]))^* (\Delta^{Wirt}_{\eta,\zeta})$$
where $T_{\tp}: P^0_\zeta \ra P'_\zeta$ denotes translation by
$\cO_{\tC}(\tp - \sigma \tp) \in P'_\zeta$. Let us choose a
square-root $\delta \in P^0_\zeta$ of 
$\cO(2\tp - 2\sigma \tp) \in P^0_\zeta$. Then we have the equality
$T^*_{\tp}\LL^\eta_\zeta = T^*_\delta \LL^\eta_\zeta$ among line
bundles over $P^0_\zeta$. We also recall the equality $[2] \circ T_\epsilon
= T_\delta \circ [2]$, where $\epsilon \in P^0_\zeta$ is a square-root
of $\delta$. Using this notation, one easily verifies that the
composite map $\phi^{-*}_{3,\zeta} \circ s^-_\eta$ factorizes as
follows (analogue of lemma \ref{lemfact3})
\begin{eqnarray} \label{lemfact3-}
H^0(P^{ev}_\eta,3\Xi_\eta)^\vee_- & \map{+T^*_\zeta \Xi_\eta} &
H^0(P^{ev}_\eta,\LL_\eta^\zeta)^\vee_-
\map{\eqref{pwd-}} H^0(P'_\zeta,\LL_\zeta^\eta)_+  \nonumber \\ 
 & \map{(T^{-1}_\delta \circ T_{\tp})^*} & H^0(P^0,\LL^\eta_\zeta)_\pm
\map{[2]^*} H^0(P^\zeta_\zeta, 8\Xi_\zeta)_\pm
\map{T^*_\epsilon} H^0(P^\zeta_\zeta, T^*_\epsilon 8\Xi_\zeta)_\pm
\end{eqnarray} 
We note that the linear isomorphism $(T^{-1}_\delta \circ T_{\tp})^*$
depends on the choice of the square-root $\delta$, which implies
the indeterminacy in the sign $\pm$ of the eigenspaces. This sign
is irrelevant for the rest of the proof. The factorization 
\eqref{lemfact3-} allows us to conclude as above.

\subsection{$n=4$}

As in the previous section, we first consider the duality on
$\MM^+(\Spin_4)$.
Let $(L,M) \in P^\zeta_\zeta \times P^\zeta_\zeta$ and  $A = 
\phi^+_{4,\zeta}(L,M) \in \MM^+(\Spin_4)$. Then we have an equality,
which is a consequence of lemma \ref{lemdec4}:
$ \forall (L,M) \in P^\zeta_\zeta \times P^\zeta_\zeta, \ 
\forall \xi \in P^{ev}_\eta \cup P^{odd}_\eta$,
$$ 
h^0(C,A(V)\otimes \pi_{\eta*} \xi) = h^0(C, \pi_{\zeta*}(LM) \otimes
\pi_{\eta*} \xi) +  h^0(C, \pi_{\zeta*}(LM^{-1}) \otimes
\pi_{\eta*} \xi) 
$$
Hence the pull-back of the divisor $\DD^+_\eta$ \eqref{defDeta} by the morphism
$$(\id \times \phi^+_{4,\zeta}) : (P^{ev}_\eta \cup P^{odd}_\eta) \times
(P^\zeta_\zeta \times P^\zeta_\zeta) \lra
(P^{ev}_\eta \cup P^{odd}_\eta) \times \MM^+(\Spin_4) $$
splits into two divisors
\begin{equation} \label{decdiveta4}
(\id \times \phi^+_{4,\zeta})^*(\DD^+_\eta) =
(\id \times m)^* \left[ p^*_{12} \Delta^{Wirt}_{\eta,\zeta}
+ p^*_{13} \Delta^{Wirt}_{\eta,\zeta} \right]
\end{equation}
where $m$ is the isogeny $P^\zeta_\zeta \times
P^\zeta_\zeta \ra P^0_\zeta \times P^0_\zeta$ defined in
\eqref{isogeny} and $p_{ij}$ is the projection of
 $(P^{ev}_\eta \cup P^{odd}_\eta) \times
 P^0_\zeta \times P^0_\zeta$ onto the $i$-th and $j$-th
factors. The decomposition \eqref{decdiveta4} leads to the
following factorization.

\begin{lem} \label{lemfact4}
For any $\eta,\zeta$ satisfying $\langle \eta,\zeta \rangle = 1$,
the linear map $s^+_\eta$ composed with the pull-back induced
by $\phi^+_{4,\zeta}$
$$H^0(P^{ev}_\eta,4\Xi_\eta)^\vee_+
\oplus H^0(P^{odd}_\eta,4\Xi_\eta)^\vee_+ \map{s^+_\eta}
H^0(\MM^+(\Spin_4), \Theta(\cc^4)) \map{\phi^{+*}_{4,\zeta}}
H^0(P^\zeta_\zeta \times P^\zeta_\zeta, 4\Xi_\zeta \boxtimes
4\Xi_\zeta)$$
factorizes as follows
\begin{eqnarray*}
 & & H^0(P^{ev}_\eta,4\Xi_\eta)^\vee_+
\oplus H^0(P^{odd}_\eta,4\Xi_\eta)^\vee_+  \\
 &  \lra & \left[ H^0(P^{ev}_\eta, \LL_\eta^\zeta)^\vee_+
\otimes H^0(P^{ev}_\eta, \LL_\eta^\zeta)^\vee_+ \right]
\oplus
 \left[ H^0(P^{odd}_\eta, \LL_\eta^\zeta)^\vee_+
\otimes H^0(P^{odd}_\eta, \LL_\eta^\zeta)^\vee_+ \right] \\
 & \map{\eqref{pwd+}} &
\left[ H^0(P^0_\zeta, \LL_\zeta^\eta)_+ \otimes
 H^0(P^0_\zeta,  \LL_\zeta^\eta)_+ \right]
\oplus
\left[ H^0(P^0_\zeta, \LL_\zeta^\eta)_- \otimes
 H^0(P^0_\zeta, \LL_\zeta^\eta)_- \right] \\
 & \map{m^*} & 
H^0(P^\zeta_\zeta \times P^\zeta_\zeta, 4\Xi_\zeta \boxtimes
4\Xi_\zeta)_+
\end{eqnarray*}
where the first map is the dual of the multiplication map.
\end{lem}

As in the $n=3$ case, it will be enough to show that
the composite map
\begin{eqnarray} \label{dualsch4}
\sum_\eta H^0(P^{ev}_\eta,4\Xi_\eta)^\vee_+
\oplus H^0(P^{odd}_\eta,4\Xi_\eta)^\vee_+
\map{\sum s^+_\eta} H^0(\MM^+(\Spin_4), \Theta(\cc^4))  \nonumber \\
 \map{\phi^{+*}_4} \sum_\zeta H^0(P^\zeta_\zeta \times
P^\zeta_\zeta, 4\Xi_\zeta \boxtimes 4\Xi_\zeta)_+ 
\end{eqnarray}
is injective. By lemma \ref{lemfact4} this linear map \eqref{dualsch4}
factorizes as follows
\begin{eqnarray*}
 & & \sum_{\eta \in JC[2]} H^0(P^{ev}_\eta,4\Xi_\eta)^\vee_+
\oplus H^0(P^{odd}_\eta,4\Xi_\eta)^\vee_+  \\
 & \map{(1)} & \sum_{\eta \in JC[2]} \sum_{\zeta \in P^0_\eta[2]}
\left[ H^0(P^{ev}_\eta, \LL_\eta^\zeta)^\vee_+
\otimes H^0(P^{ev}_\eta, \LL_\eta^\zeta)^\vee_+ \right]
\oplus
 \left[ H^0(P^{odd}_\eta, \LL_\eta^\zeta)^\vee_+
\otimes H^0(P^{odd}_\eta, \LL_\eta^\zeta)^\vee_+ \right] \\
 & \map{(2)} & \sum_{\zeta \in JC[2]} \sum_{\eta \in P^0_\zeta[2]}
\left[ H^0(P^0_\zeta, \LL_\zeta^\eta)_+ \otimes
 H^0(P^0_\zeta,  \LL_\zeta^\eta)_+ \right]
\oplus
\left[ H^0(P^0_\zeta, \LL_\zeta^\eta)_- \otimes
 H^0(P^0_\zeta, \LL_\zeta^\eta)_- \right] \\
 & \map{(3)} & \sum_{\zeta \in JC[2]} 
H^0(P^\zeta_\zeta \times P^\zeta_\zeta, 4\Xi_\zeta \boxtimes
4\Xi_\zeta)_+
\end{eqnarray*} 
Map (1) is the dual of the multiplication map in 
proposition \ref{multmap}.2, hence is injective. Map (2) is the
Prym-Wirtinger duality \eqref{pwd+}, and
map (3) is the isomorphism  \eqref{iso4+}. Hence the
composite map is injective.

\bigskip

In order to avoid repetition, we will just indicate the
changes to be done to adapt the previous proof to
$\MM^-(\Spin_4)$. The analogue of \eqref{decdiveta4}
is 
$$(\id \times \phi^-_{4,\zeta})^*(\DD^-_\eta) =
\Big(\id \times (T_{\tp} \times \id) \circ m \Big)^* 
\left[ p^*_{12} \Delta^{Wirt}_{\eta,\zeta}
+ p^*_{13} \Delta^{Wirt}_{\eta,\zeta} \right]$$
where the RHS-divisor is taken in $(P^{ev}_\eta \cup
P^{odd}_\eta) \times P'_\zeta \times P^0_\zeta$.  This
implies that the composite map $\phi^{-*}_{4,\zeta} \circ s^-_\eta$
factorizes as follows (analogue of lemma \ref{lemfact4})
\begin{eqnarray*}
 & & H^0(P^{ev}_\eta,4\Xi_\eta)^\vee_-
\oplus H^0(P^{odd}_\eta,4\Xi_\eta)^\vee_-  \\
 &  \lra & \left[ H^0(P^{ev}_\eta, \LL_\eta^\zeta)^\vee_-
\otimes H^0(P^{ev}_\eta, \LL_\eta^\zeta)^\vee_+ \right]
\oplus
 \left[ H^0(P^{odd}_\eta, \LL_\eta^\zeta)^\vee_-
\otimes H^0(P^{odd}_\eta, \LL_\eta^\zeta)^\vee_+ \right] \\
 & \map{\eqref{pwd+} \eqref{pwd-}} &
\left[ H^0(P'_\zeta, \LL_\zeta^\eta)_+ \otimes
 H^0(P^0_\zeta,  \LL_\zeta^\eta)_+ \right]
\oplus
\left[ H^0(P'_\zeta, \LL_\zeta^\eta)_- \otimes
 H^0(P^0_\zeta, \LL_\zeta^\eta)_- \right] \\
 & \map{(T^{-1}_\delta \circ T_{\tp})^* \otimes \id} &
\left[ H^0(P^0_\zeta, \LL_\zeta^\eta)_\pm \otimes
 H^0(P^0_\zeta,  \LL_\zeta^\eta)_+ \right]
\oplus
\left[ H^0(P^0_\zeta, \LL_\zeta^\eta)_\mp \otimes
 H^0(P^0_\zeta, \LL_\zeta^\eta)_- \right] \\
& \map{m^*} & 
H^0(P^\zeta_\zeta \times P^\zeta_\zeta, 4\Xi_\zeta \boxtimes
4\Xi_\zeta)_\pm \map{(T_\epsilon \times T_\epsilon)^*}
H^0(P^\zeta_\zeta \times P^\zeta_\zeta, T^*_\epsilon 4\Xi_\zeta \boxtimes
T^*_\epsilon 4\Xi_\zeta)_\pm
\end{eqnarray*}
where we used the equality $(T_\delta \times \id) \circ m = m
\circ (T_\epsilon \times T_\epsilon)$. Now we conclude as above.

\subsection{$n=5$}

Since the method of the proof is the same as for $n=3$ and
$n=4$, we will just indicate the main steps. Again we start
with $\MM^+(\Spin_5)$. As a consequence of \eqref{lemdec5}, we
see that the pull-back of the divisor $\DD_\eta^+$ by the
morphism $\id \times \phi^+_{5,\zeta}$ splits as follows
(notation as above)
\begin{equation*} \label{decdiveta5}
(\id \times \phi^+_{5,\zeta})^*(\DD^+_\eta) =
p_1^*(T^*_\zeta \Xi_\eta) + p_{23}^*((\id \times
[2]^*) (\Delta^{Wirt}_{\eta,0})) + p_{13}^*((\id \times
[2]^*) (\Delta^{Wirt}_{\eta,\zeta})) 
\end{equation*}
This decomposition implies that the composite map $\phi^{+*}_{5,\zeta}
\circ s_\eta^+$ factorizes as follows
\begin{eqnarray*}
H^0(P^{ev}_\eta,5\Xi_\eta)^\vee_+ \map{+T^*_\zeta \Xi_\eta}
H^0(P^{ev}_\eta,2\Xi_\eta)^\vee \otimes 
H^0(P^{ev}_\eta,\LL_\eta^\zeta)^\vee_+ \\
\map{\eqref{pwd+}}  
H^0(JC,2\Theta \otimes \eta)_+ \otimes 
H^0(P^0_\zeta, \LL_\zeta^\eta)_+ \map{[2]^*} 
H^0(JC,8\Theta)_+ \otimes H^0(P^\zeta_\zeta, 8\Xi_\zeta)_+
\end{eqnarray*}
As in the previous sections, the composite map
$\phi_5^{+*} \circ (\sum_{\eta} s^+_\eta)$ factorizes
\begin{eqnarray*}
\sum_{\eta \in JC[2]} H^0(P^{ev}_\eta,5\Xi_\eta)^\vee_+ & \map{(1)} &
\sum_{\eta \in JC[2]} H^0(P^{ev}_\eta,2\Xi_\eta)^\vee \otimes
H^0(P^{ev}_\eta, 3\Xi_\eta)^\vee_+ \\
 & \map{(2)} & \sum_{\eta \in JC[2]} \sum_{\zeta \in P^0_\eta[2]}
 H^0(P^{ev}_\eta,2\Xi_\eta)^\vee \otimes
 H^0(P^{ev}_\eta,\LL_\eta^\zeta)^\vee_+ \\
 & \map{(3)} & \sum_{\zeta \in JC[2]} \sum_{\eta \in P^0_\zeta[2]}
 H^0(JC, 2\Theta \otimes \eta)_+ \otimes 
 H^0(P^0_\zeta, \LL_\zeta^\eta)_+ \\
 & \map{(4)} & \sum_{\zeta \in JC[2]} 
 H^0(JC,8\Theta)_+ \otimes H^0(P^\zeta_\zeta,8\Xi_\zeta)_+ 
\end{eqnarray*}

The map (1) is the dual of the (even part) multiplication map in 
proposition \ref{multkk}.1 ($l=2,m=3$).
The map (2) is the dual of the multiplication map in proposition 
\ref{multmap}.1 tensored with $H^0(P^{ev}_\eta,2\Xi_\eta)^\vee$. Map (3) 
is Prym-Wirtinger duality \eqref{pwd+} (take $\zeta = 0$ on the
first factor),  and map (4) is an injection, by \eqref{dupl}. 
Since all four linear maps are injective, we are
done.

\bigskip
The proof of the duality for $\MM^-(\Spin_5)$ does not present
any difficulty and we leave it to the reader.

\subsection{$n=6$}

The duality for $\MM^+(\Spin_6)$ can be proved by induction
(see next section) since we know that it holds for $\MM^+(\Spin_2)$.
Since $\MM^-(\Spin_2) = \emptyset$, we have to deal with this 
case separately. Again the proof is similar to the previous
ones and we will omit it.

\subsection{the general case}

In this section we shall prove theorem \ref{duality} by induction on
$n$. We will also denote by $\Theta(V)$ the pull-back by $\phi^\pm_n$
of the determinant line bundle $\Theta(V) = \Theta(\cc^n)$ 
to the Schottky variety $SCH_n$. Our induction hypothesis $\HH_n$ will
be the following 

\bigskip
\noindent
$\HH_n$ for $n$ odd:
$$\phi^{\pm*}_n \circ (\sum s^\pm_\eta): \  \sum_\eta
H^0(P^{ev}_\eta,n\Xi_\eta)^\vee_\pm \lra H^0(SCH_n,\Theta(V)) \ \   
\text{is injective}$$

\bigskip
\noindent
$\HH_n$ for $n$ even:
$$\phi^{\pm*}_n \circ (\sum s^\pm_\eta): \  \sum_\eta
H^0(P^{ev}_\eta,n\Xi_\eta)^\vee_\pm \oplus 
H^0(P^{odd}_\eta,n\Xi_\eta)^\vee_\pm \lra H^0(SCH_n,\Theta(V)) \ \ 
\text{is injective}$$

Since (\cite{oxbury2} theorem 3.1) the LHS-space and 
$H^0(\MM^\pm(\Spin_n),\Theta(V))$ have the same dimension (see remark 3
of section 8 for $\MM^-(\Spin_n)$),
the assumption $\HH_n$ implies theorem 1.1. We already
proved $\HH_n$ for $2 \leq n \leq 6$. Let us assume $\HH_n$ and prove
$\HH_{n+4}$. First we easily verify that the pull-back of the
divisor $\DD^\pm_{n+4,\eta}$ under the natural map (e.g. for $n$ odd)
(see lemma 3.1)
$$ \id \times \iota: \ P^{ev}_\eta \times \MM^+(\Spin_4) \times
\MM^\pm(\Spin_n) \lra P^{ev}_\eta \times \MM^\pm(\Spin_{n+4}) $$
splits into two divisors
\begin{equation} \label{decDetan}
(\id \times \iota)^* (\DD^\pm_{n+4,\eta}) = p^*_{12}(\DD^+_{4,\eta}) +
p^*_{13}(\DD^\pm_{n,\eta}).
\end{equation}
We distinguish two cases.

\subsubsection{$n$ odd}

As a consequence of the decomposition \eqref{decDetan} and the proof
of the $n=4$ case, the map $\phi^{\pm*}_{n+4} \circ (\sum s^\pm_\eta)$
factorizes as follows
\begin{eqnarray*}
 & & \sum_{\eta \in JC[2]} H^0(P^{ev}_\eta,(n+4)\Xi_\eta)^\vee_\pm \\
 & \map{(1)} &
\sum_{\eta \in JC[2]} H^0(P^{ev}_\eta,4\Xi_\eta)^\vee_+ \otimes
H^0(P^{ev}_\eta, n\Xi_\eta)^\vee_\pm \\
 & \map{(2)} & \sum_{\eta \in JC[2]} \sum_{\zeta \in P^0_\eta[2]}
 H^0(P^{ev}_\eta,\LL_\eta^\zeta)^\vee_+ \otimes
 H^0(P^{ev}_\eta,\LL_\eta^\zeta)^\vee_+ \otimes
 H^0(P^{ev}_\eta, n\Xi_\eta)^\vee_\pm  \\
 & \map{(3)} & \sum_{\zeta \in JC[2]} \sum_{\eta \in P^0_\zeta[2]}
 H^0(P^0_\zeta, \LL_\zeta^\eta)_+ \otimes
  H^0(P^0_\zeta, \LL_\zeta^\eta)_+ \otimes
 H^0(P^{ev}_\eta, n\Xi_\eta)^\vee_\pm  \\
  & \map{(4)} & \sum_{\zeta \in JC[2]} \left[ \sum_{\eta \in P^0_\zeta[2]}
 H^0(P^0_\zeta, \LL_\zeta^\eta)_+ \otimes
  H^0(P^0_\zeta, \LL_\zeta^\eta)_+ \right] \otimes
 \left[ \sum_{\eta \in JC[2]} H^0(P^{ev}_\eta, n\Xi_\eta)^\vee_\pm
 \right]   \\
 & \map{(5)} & \sum_{\zeta \in JC[2]} 
 H^0(P^\zeta_\zeta \times P^\zeta_\zeta, 4\Xi_\zeta \boxtimes
 4\Xi_\zeta) \otimes H^0(SCH_n, \Theta(V)) = H^0(SCH_{n+4},\Theta(V)) 
\end{eqnarray*}
The maps are as follows: (1) is the dual of the multiplication map
in proposition \ref{multkk}.2 ($l=2,m=n$). If $n=2$, we use 
proposition \ref{multkk}.1 ($l=2,m=4$). (2) is the dual of 
the (even part) multiplication map in proposition \ref{multmap}.2. (3)
is the Prym-Wirtinger duality \eqref{pwd+} tensored with
$H^0(P^{ev}_\eta,n\Xi_\eta)^\vee_\pm$. (4) is the inclusion
$H^0(P^{ev}_\eta,n\Xi_\eta)^\vee_\pm \inj \sum_\eta 
H^0(P^{ev}_\eta,n\Xi_\eta)^\vee_\pm$.
(5) is an injection coming from the direct sum decomposition 
\eqref{iso4+} tensored with the injective map of $\HH_n$. 
Finally the last equality follows from the definition of 
$SCH_{n+4}$. Since all linear maps are injective, the
composite map is injective and we have proved $\HH_{n+4}$.

\subsubsection{$n$ even}

Since this case is similar to the odd case, we just indicate the
minor changes to be done on the sequence of maps (1),...,(5).
We take into account the additional factors $\sum_\eta
H^0(P^{odd},(n+4)\Xi_\eta)^\vee_\pm$, for which we can write down
the maps (1),...,(4) with $P^{ev}_\eta$ replaced by $P^{odd}_\eta$
and $H^0(P^0_\zeta,\LL^\eta_\zeta)_+$ by 
$H^0(P^0_\zeta,\LL^\eta_\zeta)_-$ (see \eqref{pwd+}). Thus
adding the two copies of map (5) (written for $P^{ev}$ and $P^{odd}$),
we observe that we still have an injection coming from 
\eqref{iso4+}.

\section{Final remarks}

{\bf 1.} Let us briefly indicate why we used the name ``Schottky variety''
for the products of Prym varieties. Indeed, for $n=3$, the image 
of $\phi_3$ consists of a union of Kummer varieties of Pryms (resp. 
of the Jacobian) which intersect along some $4$-torsion points (resp. 
$2$-torsion points). We refer to \cite{VGP1}, \cite{donagi1},
\cite{pauly} for a proof of these intersection properties, which may
also be deduced from lemma \ref{lemdec3}. The coordinates of the
intersection points may be interpreted as the famous Schottky-Jung
identities among theta-constants. Let us call the image
$\phi_3(SCH_3) \subset \mod0$ the Schottky configuration, which we
embed in projective space $\pp^{2^g-1} = |\LL|^\vee$, where $\LL$
is the ample generator of $\pic(\mod0)$. In particular, we
have $\Theta(\cc^3) = \LL^4$. Then we can deduce from the
injectivity of \eqref{dualsch3} 

\begin{cor}
If a quartic in $\pp^{2^g-1}$ vanishes on the Schottky configuration,
then it vanishes on the whole of $\mod0$.
\end{cor}

This corollary was already proved in \cite{VGP1} corollary 2. Moreover
we have been able to show that $\mod0 \subset \pp^{2^g-1}$ is defined
by quartic equations. One can thus recover $\mod0$ from the
Schottky configuration purely geometrically. One might speculate
whether, given the Schottky configuration of abelian varieties,
one might reconstruct the curve $C$. This is the ``small Schottky''
conjecture, see \cite{donagi2}.

\bigskip
\noindent
{\bf 2.} As observed in \cite{oxbury2} remark 1.2 and \cite{ls} section
7.10,
the reduced divisors 
$D^{(n)}_\kappa$, with  $\kappa \in \vartheta(C)$,
whose  support is given by
$$ \supp \ D^{(n)}_\kappa = \{ A \in \MM^+(\Spin_n) \ | \ h^0(C, A(V) \otimes 
\kappa ) >0 \}, $$
are elements of the linear system $\pp H^0(\MM^+(\Spin_n),\mathcal{P})$.
Here $\mathcal{P}$ is a Pfaffian square-root of $\Theta(V)$,i.e.,
$\mathcal{P}^2 = \Theta(V)$. To avoid existence problems, we work
over the moduli stack. Then we shall prove

\begin{prop}
A basis of the linear system $|\mathcal{P}|$ is given by the divisors
$D_\kappa^{(n)}$ with $\kappa \in \vartheta(C)$, if $n$ is even, and
$\kappa \in \vartheta(C)$ with $\epsilon(\kappa) =1$, if $n$ is odd.
\end{prop}

\begin{proof}
It is enough to prove linear independence, which is done by
induction on $n$. For the first cases we refer to \cite{beauville}
proposition A.8 ($n=2$), theorem 1.2 ($n=3$), proposition A.5 ($n=4$).
The statement for $n=5$ will follow by pulling back $D_\kappa^{(5)}
\subset \MM^+(\Spin_5)$ under the map $\MM^+(\Spin_3) \times
\MM^+(\Spin_2) \map{\iota} \MM^+(\Spin_5)$. We observe that
$\iota^*(D_{\kappa}^{(5)}) = p_1^*(D^{(3)}_\kappa) +
p_2^*(D^{(2)}_\kappa)$. Since the $D^{(3)}_\kappa$ form a basis,
we can conclude. Let us assume that the $D^{(n)}_\kappa$'s are
independent. Then the equality $\iota^*(D^{(n+4)}_\kappa) =
p_1^*(D^{(n)}_\kappa) + p_2^*(D^{(4)}_\kappa)$ implies that
the $D^{(n+4)}_\kappa$'s are also independent.
\end{proof}
 
\bigskip
\noindent
{\bf 3.} The ``twisted'' Verlinde formula given in \cite{oxbury2}
conjecture 1.1, which computes the dimension of 
$H^0(\MM^-(\Spin_n),\Theta(V))$,
can be deduced from a forthcoming work by Y. Laszlo and C. Sorger.
Following the techniques of the paper \cite{ls}, the authors show that
the ``twisted'' Verlinde space can be identified to a conformal
block, whose dimension (worked out by C. Woodward) is given by
Oxbury's formula. 

\bigskip
\noindent
{\bf 4.} We can prove a refinement of theorem 1.1 by considering 
the linear action of the Heisenberg group $\Heis$ on
$H:= H^0(\MM^\pm(\Spin_n),\Theta(V))$. By induction, we check
that $H$ is a $\Heis$-module of level $4$ for the $JC[2]$-action
described in \eqref{galaction}. Hence the linear action 
of $\Heis$ factors through $JC[2]$ and we can consider its
character space decomposition $H = \sum_{\eta \in JC[2]} H_\eta$.
One proves that the image of $s^\pm_\eta$ is contained in
$H_\eta$. Since $\sum_\eta s_\eta$ is surjective, we get
equality $\im(s_\eta) = H_\eta$. So, we have isomorphisms, e.g. for
$n$ odd,
$$ \forall \eta \in JC[2] \qquad
s_\eta^\pm \ : \ H^0(P^{ev}_\eta, (2m+1)\Xi_\eta)^\vee_\pm
\map{\sim} H^0(\MM^\pm(\Spin_{2m+1}), \Theta(\cc^{2m+1}))_\eta $$

\flushleft{C. Pauly \\
Laboratoire J.-A. Dieudonn\'e \\
Universit\'e de Nice Sophia Antipolis \\
Parc Valrose \\
06108 Nice Cedex 02 \\ France \\
email: pauly@math.unice.fr}

\bigskip
\noindent
\flushleft{S. Ramanan \\
Tata Institute of Fundamental Research \\
Homi Bhabha Road \\
Bombay 400005 \\
India\\
email: ramanan@math.tifr.res.in}


\begin{thebibliography}{999999}
\bibitem[B1]{beau1} A. Beauville: Fibr\'es de rang $2$ sur une 
courbe, fibr\'es d\'eterminant et fonctions th\^eta, Bull. Soc.
Math. France 116 (1988), 431-448

\bibitem[B2]{beauville} A. Beauville: Fibr\'es de rang $2$ sur une
courbe, fibr\'es d\'eterminant et fonctions th\^eta II, Bull. Soc.
Math. France 119 (1991), 259-291

\bibitem[B3]{beau3} A. Beauville: Vector bundles on curves and
generalized theta functions: recent results and open problems,
Current topics in complex algebraic geometry (Berkeley, CA, 1992/93),
Math. Sci. Res. Inst. Publ. 28, 17-33

\bibitem[BNR]{bnr} A. Beauville, M.S. Narasimhan, S. Ramanan:
Spectral curves and the generalized theta divisor, J. reine angew.
Math. 398 (1989), 169-179 

\bibitem[D1]{donagi1} R. Donagi: Non-Jacobians in the Schottky loci,
Ann. of Math. 126 (1987), 193-217

\bibitem[D2]{donagi2} R. Donagi: Big Schottky, Invent. Math. 89 (1987),
569-599

\bibitem[DT]{dt} R. Donagi, L. Tu: Theta functions for $\SL(n)$ versus
$\GL(n)$, Math. Res. Lett. 1, no.3 (1994), 345-357 


\bibitem[vGP]{VGP1} B. van Geemen, E. Previato: Prym varieties and
the Verlinde formula, Math. Ann. 294 (1992), 741-754


\bibitem[K]{kempf} G. Kempf: Multiplication over abelian varieties,
Amer. J. Math. 110 (1988), 765-773

\bibitem[Kh]{khaled} A. Khaled: Equations des vari\'et\'es de Kummer,
Math. Ann. 295 (1993), 685-701

\bibitem[M1]{mum1} D. Mumford: Prym varieties I, In Contributions
to Analysis, Academic Press, New York (1974)

\bibitem[M2]{mum2} D. Mumford: Theta-characteristics of an algebraic
curve, Ann. Sci. Ec. Norm. Sup. 4 (1971), 181-192

\bibitem[M3]{mum3} D. Mumford: On the equations defining abelian varieties,
Invent. Math. 1 (1966), 287-354 

\bibitem[LS]{ls} Y. Laszlo, C. Sorger: The line bundle on the
moduli of parabolic $G$-bundles over curves and their sections,
Ann. Sci. Ec. Norm. Sup. 30 (1997), 499-525

\bibitem[O1]{oxbury1} W.M. Oxbury: Prym varieties and the moduli
of spin bundles, Algebraic geometry, ed. P.E. Newstead, Lect. Notes
in Pure App. Math. vol. 200 (Marcel Dekker 1998), 351-376

\bibitem[O2]{oxbury2} W.M. Oxbury: Spin Verlinde spaces and Prym
theta functions, Proc. London Math. Soc. (3) 78 (1999), 52-76

\bibitem[OP]{op} W.M. Oxbury, C. Pauly: $SU(2)$-Verlinde
spaces as theta spaces on Pryms, Int. J. Math. 7 (3) (1996), 
393-410

\bibitem[OW]{ow} W.M. Oxbury, S.M.J. Wilson: Reciprocity laws
in the Verlinde formulae for the classical groups, Trans. Amer. Math.
Soc. 348 (1996), 2689-2710

\bibitem[P]{pauly} C. Pauly: On Pryms, rank $2$ bundles and nonabelian
theta functions, Int. J. Math. 8 (2) (1997), 267-287

\bibitem[R1]{ramanan} S. Ramanan: Orthogonal and spin bundles over 
hyperelliptic curves, in Geometry and Analysis, papers dedicated to 
V.K. Patodi, Springer 1981

\bibitem[R2]{ram2} S. Ramanan: The Schottky configuration, preprint 1995

\bibitem[S]{sorger} C. Sorger: La formule de Verlinde, S\'eminaire 
Bourbaki, Vol. 1994/5, Ast\'erisque 237 (1996), Exp. No. 794, 3, 87-114

\end{thebibliography}
\end{document}